\documentclass{article}

\usepackage{amssymb,dsfont,latexsym}

\newtheorem{thm}{Theorem}

\newtheorem{prop}[thm]{Proposition}
\newtheorem{lem}[thm]{Lemma}
\newtheorem{cor}[thm]{Corollary}

\newcommand\enu[1]{\smallskip\newline\makebox[5mm][l]{\rm(#1)}}

\newcommand\bp{\noindent{\it Proof.}\ }

\newcommand\Ad{{\rm Ad}\,}

\begin{document}

\author{Erling St{\o}rmer}

\title{Tensor products of positive maps of matrix algebras}

\maketitle

\begin{abstract}
We give conditions for when the tensor product of two positive maps
between matrix algebras is a positive map.  This happens when one map
belongs to a symmetric mapping cone and the other to the dual cone.
Necessary and sufficient conditions are given.
\end{abstract}

\section*{Introduction}
Even though positive linear maps appear in many situations in
operator algebras and quantum information theory, main attention has
so far been on completely positive maps.  One reason for this is that
tensor products of completely positive maps are positive, while this
is false for general positive maps. In the present paper we shall
consider the problem of when the tensor product of two positive maps
is positive in the case when the underlying Hilbert spaces are finite
dimensional. 

It turns out that the problem is intimately related to symmetric
mapping cones and their dual cones, see below for definitions.  More
specifically, if $\1C$ is a symmetric mapping cone in $P(H)$ - the
positive linear maps of the bounded operators $B(H)$ on $H$ into
itself - then a map $\phi$ belongs to the dual cone $\1C^o$  of $\1C$ if if
and only if $\psi\otimes\phi$ is positive for all $\psi\in \1C$.
Indeed, it suffices to know that $\psi\otimes\phi(p)\geq 0$, where
$\frac{1}{n} p, n=dim H$ is the density matrix for the maximally entangled
state.  As an application we show that if $K$ and $L$ are two other
finite dimensional Hilbert spaces, and $\psi\colon B(K) \to B(H)$,
$\phi\colon B(L) \to B(H)$, then $\psi\otimes\phi $ is positive when
$\psi$ is $\1C$-positive, and $\phi$ is  $\1C^o$-positive.

We now recall the main concepts encountered in the sequel.  By a
\textit{mapping cone} $\1C$ we mean a closed subcone of $P(H)$ such
that if $\phi\in \1C$ and $\alpha, \beta  \in CP(H)$- the completely
positive maps in $P(H)$, then $\alpha\circ\phi\in \1C$ and
$\phi\circ\beta\in \1C$.  We say $\1C$ is \textit{symmetric} if
$\phi\in\1C$ if and only if $\phi^*\in\1C$ if and only if $\phi^t =
t\circ \phi\circ t \in\1C$, where $\phi^*$ is the adjoint map of
$\phi$ in the Hilbert-Schmidt structure on $B(H),$ viz. $Tr(\phi(a)b) =
Tr(a\phi^*(b))$ for $a,b\in B(H)$, and $t$ is the transpose map on
$B(H)$ with respect to an orthonormal basis $e_1,...,e_n$ for $H$, and
$Tr$ is the usual trace on $B(H)$.

If $\phi\colon B(K)\to B(H) = B(H\otimes H)$ then the functional $\widetilde \phi $ on
$B(K)\otimes B(H)$ defined by 
$$
\widetilde \phi(a\otimes b) = Tr(\phi(a) b^t),
$$
plays an important role in the theory.  For example, $\widetilde \phi$
is positive if and only if $\phi$ is completely positive \cite{s1}
.  By \cite{s2} and \cite{s3} its
density matrix is the transpose $\1C_{\phi}^t$ of the Choi matrix $\1C_{\phi}$ for
$\phi$, defined by 
$$
\1C_{\phi}   = \sum_{i,j =1} ^{n}  e_{ij} \otimes \phi(e_{ij} ) =
  \iota\otimes \phi (p),
$$
where $\iota$ is the identity map and  $p=\sum_{ij} e_{ij} \otimes e_{ij},$ and $(e_{ij})$ is a
complete set of matrix units for $B(H)$ such that $e_{ij} e_k =
\delta_{jk} e_i$, see \cite{C} .

Let $\1C$ be a mapping cone in $P(H)$. Then its dual cone is defined
by
$$
\1C^o = \{\phi\in P(H): Tr(C_{\phi} C_{\psi}) \geq 0,  \forall
\psi\in \1C \}.
$$
If $\1C$ is symmetric, then by  \cite{s4}            $\1C^o$ is also a
symmetric mapping cone.  We refer to the books  \cite{ER} and \cite{VP}
for the theory of completely positive maps.

Most of this work was done during a visit to Institute Mittag-Leffler
(Djursholm, Sweden)

\section{The main results}

   Let $\pi\colon B(H)\to B(H)$ be defined by $\pi(a\otimes b) = b^t
   a$.  By \cite{s3}     Lemma 10, if $\phi\in P(H)$ then 
\begin{equation} \label{e1}
\widetilde \phi = Tr \circ \pi\circ (\iota\otimes \phi^{*t})   
\end{equation}  
In particular
\begin{equation}\label{e2}
\widetilde \iota(x) = Tr\circ \pi (x) = Tr(C_{\iota} x) = Tr(p x).
\end{equation}
Thus
\begin{eqnarray*}
\widetilde\phi(x)
 &=& Tr\circ \pi(\iota\otimes \phi^{*t}(x))\\
&=& Tr( p (\iota\otimes \phi^{*t})(x)) \\
&=& Tr(\iota\otimes \phi^t(p) x).
\end{eqnarray*}

\begin{lem}\label{lem1}
Let $\psi,\phi \colon B(H)\to B(H)$.  Then
\enu{i} $ ( \phi\circ\psi)\widetilde   (x) = Tr(\psi^* \otimes \phi^t (p)
x), \forall  x\in B(H\otimes H)$.
\enu{ii} $\psi^{*t} \otimes \phi (p) = \iota \otimes
(\phi\circ\psi)(p)$.
\end{lem}

\bp
By equations (1) and (2)
\begin{eqnarray*}
(\phi\circ\psi)\widetilde  (a\otimes  b)
&=& Tr\circ\pi (\iota\otimes (\phi\circ\psi)^{*t} (a\otimes b))\\
&=& Tr\circ\pi(a\otimes (\psi^* \circ \phi^* (b^t)) ^t )\\
&=& Tr(a (\psi^* \circ \phi^*) (b^t))\\
&=& Tr(\psi(a) \phi^* (b^t))\\
&=& Tr\circ\pi(\psi(a)\otimes \phi^{*t}(b))\\
&=& Tr(p \psi(a) \otimes \phi^{*t}(b))\\
&=& Tr(\psi^* \otimes\phi^t(p)( a\otimes b)),
\end{eqnarray*}
proving (i).
Using equation (1) we also have
$$
(\phi\circ\psi)\widetilde  (x) 
= Tr(p(\iota\otimes (\phi\circ\psi)^{*t} (x))
= Tr(\iota\otimes(\phi\circ\psi)^t(p) x).
$$
Thus by (i) we have
$$
\psi^*\otimes \phi^t(p) = \iota\otimes \phi^t \circ\psi^t(p).
$$
Since this holds for all $\phi$ and $\psi$, it also holds for all
$\phi^t$ and $\psi^*$.  Thus
$$
\psi^{*t} \otimes \phi (p) = \iota \otimes \phi\circ\psi (p),
$$
completing the proof of the lemma.

We can now prove our main result.

\begin{thm}\label{thm2}
Let $\phi\in P(H)$.  Let $\1C$ be a symmetric mapping cone in $P(H)$.
Then the following conditions are equivalent.
\enu{i} $\phi\in \1C^o$  - the dual cone of $\1C$,
\enu{ii} $\phi\circ\psi$ is completely positive for all $\psi\in \1C$,
\enu{iii} $\psi\otimes\phi$ is positive  for all $\psi\in \1C$,
\enu{iv} $\psi\otimes\phi(p) \geq 0$  for all $\psi\in \1C$.
\end{thm}

\bp
Clearly (iii) implies (iv).  Since $\phi\circ\psi$ is completely 
positive if and only if $\iota\otimes \phi\circ\psi(p)\geq 0,$ by
Lemma 1  $\phi\circ\psi$ is completely positive
 if and only if 
$\psi^{*t}\otimes \phi(p)\geq 0$.  Since $\1C$ is symmetric,
$\psi\in\1C$ if and only if $\psi^{*t} \in\1C$.  Hence (ii)
is equivalent to (iv). By \cite{s4} Thm.2  a map belongs to $\1C$ if
and only if it is $\1C$-positive.  Hence by  \cite{s3}    Thm.1,  $\phi\in\1C^o$
if and only if $\psi^t\circ \phi$ is completely positive for all
$\psi\in\1C$. Since $\1C$ is symmetric this holds if and only if  
$\psi\circ\phi$ is completely positive for all $\psi\in\1C$, hence if
and only if $\phi^*\circ\psi^* = (\psi\circ\phi)^*$ is completely
positive for all $\psi^*\in \1C$, hence if and only if
$\phi^*\circ\psi$ is completely positive for all $\psi\in\1C$. But
$\1C^o$ is symmetric by  \cite{s4}, Thm.1, so $\phi\in\1C^o$ if and only
if $\phi^*\in\1C^o$.  Thus (i)$ \Leftrightarrow $(ii). 

It remains to show (i) implies (iii).  For this let $(e_i)$ be an
orthonormal basis for $H$ such that $e_{ij}e_j = e_i$, so $\frac{1}{n}
p$ is the projection onto the subspace spanned by $\sum_i e_i\otimes
e_i$.  Let $x\in H\otimes H$  Then $x=\sum_i e_i \otimes x_i$ with
$x_i\in H$.  Then there is $v\in H$ such that $ve_i = x_i$, hence
$1\otimes v (\sum_i e_i\otimes e_i) = x.$ Let $q$ be the projection
onto $\0C x$.  Then it follows that $Ad(1\otimes v)( p) = \lambda q$ for
some $\lambda >0$. 

 We have just shown that given a 1-dimensional projection $q\in B(H)$
 there exists $v \in H$ such that 
$$
1\otimes Adv (\frac{1}{n} p) = q.
$$
assuming (i) $\phi\circ\Ad v \in \1C^o$, since $ \1C^o$ is a mapping
cone.  By Lemma 1
$$
\psi^{*t} \otimes (\phi\circ Ad v) (p) = \iota\otimes \phi(Ad v \circ
\psi)(p).
$$
Since $Adv \circ \psi \in \1C$, by the equivalence of (i) and (ii)
$\phi\circ Ad v \circ\psi$ is completely positive, hence
$$
\iota\otimes \phi\circ Ad v \circ\psi(p)\geq 0.
$$
Thus by the choice of $v$, $\psi^{*t} \otimes \phi (q)\geq 0$.  Since
$q$ is an arbitrary 1-dimensional projection, $\psi^{*t} \otimes \phi$
is positive.  Again, since $\1C$ is symmetric, $\psi\circ\phi$ is
positive for all $\psi\in\1C.$  Thus (i) implies (iii), and the proof
is complete.

\medskip 
Recall that a map $\phi\colon B(K)\to B(H)$ is \textit {\1C-positive}
for a mapping cone $\1C$ if the functional $\widetilde\phi$ is
positive on the cone $ \{x\in B(K\otimes H): \iota\otimes\alpha(x) \geq
0, \forall \alpha\in \1C\}$. By \cite{s4}, Thm. 2 or \cite{s1},
Thm. 3.6, this is equivalent to $\phi$ belonging to the cone generated
by maps of the form $\alpha\circ\beta$ with $\alpha\in \1C$ and
$\beta\colon B(K)\to B(H)$ completely positive.  Recall from \cite{s4},
Thm. 1, that if $\1C$ is symmetric, so is $\1C^o$. Using these facts we
can extend the implication (i) $\Rightarrow$ (iii) in Theorem 2 to the
following more general case. 

\begin{cor}\label{cor3}
Let $K,L$ and $H$ be finite dimensional Hilbert spaces.  Let $\1C$ be
a symmetric mapping cone in $P(H)$.  Suppose $\psi\colon B(K)\to B(H)$
is $\1C$-positive, and  $\phi\colon B(L)\to B(H)$ is
$\1C^o$-positive.  Then $\psi\otimes\phi\colon B(K\otimes L)\to
B(H\otimes H)$ is positive.
\end{cor}
\bp
By the above discussion it suffices to show the corollary for $\psi$
and $\phi$ of the form $\psi=\alpha\circ\beta, \alpha\in \1C,
\beta\colon B(K)\to B(H)$ completely positive, and
$\phi=\gamma\circ\delta$ with $\gamma\in\1C^o, \delta\colon B(L)\to
B(H)$ completely positive.  Thus 
$$
\psi\otimes\phi = (\alpha\otimes\gamma) \circ (\beta\otimes \delta),
$$
is positive, since $\beta\otimes \delta$ is completely positive and
$\alpha\otimes\gamma$ is positive by Theorem 2. The proof is complete.
\medskip

\textit{Remark}.  If $\psi\colon B(K_1)\to B(H_1)$ is k-positive, i.e. $\psi\in
P_k$ in the notation of \cite{SSZ}, and $\phi\colon B(K_2)\to B(H_2)$ is
k-superpositive,  i.e. $\phi\in SP_k$ is of the form $\sum_i Ad V_i,
  V_i\colon K_2\to H_2$, then they remain the same as maps into $B(H)$
  if $H$ is a Hilbert space containing $H_1$ and $H_2$ as subspaces.
  Since $P_k^o = SP_k$, see e.g.\cite{SSZ}, it follows from Corollary
  3 that $\psi\otimes\phi$ is positive.

\medskip

In Theorem 2 it is sometimes enough to consider only one map $\psi\in
\1C$ to conclude that $\phi\in \1C^o$.  The next corollary is of this
type.

\begin{cor}\label{cor4}
Let $\psi\in P(H)$ satisfy $\psi = \psi^* = \psi^t$.  Let $\1C$ denote
the mappng cone generated by $\psi$.  Let $\phi\in P(H)$.  Then
$\phi\in \1C^o$ if and only if $\psi \otimes\phi$ is positive.
\end{cor}

\bp
$\1C$ is generated as a cone by maps of the form $Ad u\circ\psi\circ
Ad  v$, so the assumptions on $\psi$ imply that $\1C$ is a symmetric
mapping cone.  Since 
$$
Ad u\circ\psi\circ Ad v\otimes \phi = (Ad u\otimes
\iota)\circ(\psi\otimes\phi)\circ(Ad v \otimes \iota),
$$
and $ Ad u\otimes\iota$ and $Ad v \otimes \iota$ are positive maps, it
follows that $\alpha\otimes \phi$ is positive for all $\alpha\in \1C$
if and only if  $\psi \otimes\phi$ is positive, hence by Theorem 2,
$\phi\in \1C^o$ if and only if $\psi \otimes\phi$ is positive, proving
the corollary.
\medskip

\textit{Remark}.  Theorem 2 and Corollary 4 contain  well known
characterizations of completely maps.  It $\psi = \iota$ then it
satisfies the conditions of Corollary 4, so the mapping cone $\1C$
generated by $\psi$ is the cone of completely positive maps. Hence if
$\phi\in P(H)$, then by Corollary 4, $\phi\in \1C^o$ if and only if
$\iota\otimes\phi$ is positive, if and only if $\phi\in\1C$ by
definition of $\1C$, so $ \1C^o =\1C$.  By Theorem 2 we have
\begin{eqnarray*}
\1C_{\phi} = \iota\otimes\phi (p)\geq 0 
&\Leftrightarrow&(\alpha\otimes\iota)(\iota\otimes\phi (p))\geq 0,
\forall \alpha\in\1C \\
&\Leftrightarrow & \alpha\otimes\phi(p)\geq 0,  \forall \alpha\in\1C \\
&\Leftrightarrow & \phi\in \1C^o = \1C.
\end{eqnarray*}

In Corollary 4 we assumed $\psi = \psi^* = \psi^t$.  These conditions
can be easily verified by checking the corresponding conditions for
the Choi matrix.

\begin{prop}\label{prop5}
Let $\phi\in P(H)$.  Then $\phi=\phi^*=\phi^t$ if and only if
$C_{\phi}$ is a real symmetric matrix invariant under the flip
$a\otimes b \to b\otimes a$ on $B(H)\otimes B(H)$.
\end{prop}

\bp
Let $J$ be the conjugation on $H\otimes H$ defined by 
$$
Jze_i\otimes e_j = \overline{z} e_j\otimes e_i, z\in \0C,
$$
where $e_i,...,e_n$ is an orthonormal basis such that $e_{ij} e_k =
\delta_{jk} e_i$.  Then an easy computation shows that if $a,b$ are
real matrices in $B(H)$, then $Ja\otimes bJ =  b\otimes a$, so for
  $x\in B(H\otimes H)$, a real matrix with respect to the basis
  $(e_i\otimes e_j)$ for $H\otimes H$,  then $x\to JxJ$ is the flip $F$ applied to $x$.

We have $C_{\phi^t} = C_{\phi}^t,$ so $\phi=\phi^t$ if and only if
$C_{\phi} =C_{\phi}^t,$ i.e. $C_{\phi}$ is symmetric.  Since $\phi\in
P(H)$, $C_{\phi}$ is self-adjoint, hence $C_{\phi}$ is symmetric if
and only if $C_{\phi}$ is real symmetric.  Hence $\phi=\phi^t$ if and
only if $C_{\phi}$ is real symmetric. By \cite{s4}, Lem. 3,
$C_{\phi^*} =JC_{\phi}J$.  Hence $\phi=\phi^* =\phi^t$ if and only if 
 $C_{\phi}$ is real symmetric, and by the first part of the proof,
 $C_{\phi} =F(C_{\phi})$, so invariant under the flip, completing the
 proof.
\medskip.

\textit{Example}.  A specific example of a map as in Proposition 5 is given by
$\phi = Ad V$, where $V$ is a real symmetric matrix.  Indeed, for
general $V$ we have the formulas:
$$
(Ad V)^* = Ad V^*,  (Ad V)^t = Ad \overline{V},
$$
where $\overline{V} =(\overline{a_{ij}})$ if $V=(a_{ij}).$  Thus, if
$V$ is real symmetric, then $Ad V =(Ad V)^*=(Ad V)^t$.  Furthermore,
if $V$ is real symmetric and $F$ the flip then
$$
C_{AdV} = \sum_{kl} e_{kl}\otimes AdV(e_{kl}) = \sum_{ijkl} v_{ki}
v_{lj} e_{kl}\otimes e_{ij}.
$$
Thus
\begin{eqnarray*}
F(C_{AdV} )&=& F( \sum_{ijkl} v_{ki} v_{lj} e_{kl}\otimes e_{ij})\\
&=& \sum_{ijkl} v_{ki} v_{lj} e_{ij}\otimes e_{kl}\\
&=& \sum_{ijkl} v_{ik} v_{jl} e_{kl}\otimes e_{ij}\\
&=& \sum_{ijkl} v_{ki} v_{lj} e_{kl}\otimes e_{ij}\\
&=& C_{AdV} ,
\end{eqnarray*}
where we at the third equality sign changed the roles of i and k, and
l and j, and used that $V$ was symmetric at the fourth. It follows
that $C_{AdV} $ is invariant under the flip.

Department of Mathematics, University of Oslo, 0316 Oslo, Norway.

e-mail erlings@math.uio.no

\end{document}